\documentclass{article}
\usepackage{graphicx}
\usepackage{amsthm,amsfonts}
\usepackage{amsmath}
\usepackage{hyperref}

\newtheorem{thm}{Theorem}[section]
\newtheorem{theorem}[thm]{Theorem}

\newtheorem{lem}[thm]{Lemma}
\newtheorem{prop}[thm]{Proposition}
\theoremstyle{definition}

\newtheorem{rem}[thm]{Remark}
\newtheorem{defn}[thm]{Definition}
\newtheorem*{ack}{Acknowledgement}

\def \build#1#2#3{\mathrel{\mathop{\kern 0pt#1}\limits_{#2}^{#3}}}

\newcommand{\refT}[1]{Theorem~\ref{#1}}

\newcommand{\refL}[1]{Lemma~\ref{#1}}

\newcommand{\refS}[1]{Section~\ref{#1}}

\newenvironment{romenumerate}{\begin{enumerate}
 }{\end{enumerate}}

\newcounter{thmenumerate}

\newcommand\bbR{\mathbb R}

\newcommand\E{\operatorname{\mathbb E{}}}
\renewcommand\P{\operatorname{\mathbb P{}}}

\newcommand\bigpar[1]{\bigl(#1\bigr)}
\newcommand\Bigpar[1]{\Bigl(#1\Bigr)}

\newcommand\smallpar[1]{(#1)}
\def\rompar(#1){\textup(#1\textup)}    
\newcommand\xfrac[2]{#1/#2}
\newcommand\parfrac[2]{\Bigpar{\frac{#1}{#2}}}

\newcommand\ceil[1]{\lceil#1\rceil}
\newcommand\floor[1]{\lfloor#1\rfloor}

\def \ISE {\mathcal J}

\newcommand{\pa}[1]{\left(#1\right)}
\newcommand{\cro}[1]{\left[#1\right]}

\newcommand{\cov}[1]{\textnormal{Cov}\left(#1\right)}

\newcommand{\esp}[1]{\mathbb{E}\left[#1\right]}
\newcommand{\pr}[1]{\mathbb{P}\left(#1\right)}

\newcommand{\loi}{\build{=}{}{law}}

\title{The center of mass of the ISE and the Wiener index of trees}

\date{}
\author{
Svante Janson\thanks{Departement of Mathematics,
                           Uppsala University,
                           P. O. Box 480,
                           S-751 06 Uppsala, Sweden.
                       Email: {\tt svante.janson@math.uu.se}}
\and Philippe Chassaing\thanks{Institut \'{E}lie Cartan,
                                Universit\'e Henri Poincar\'e-Nancy I,
                                BP 239, 54506
                                Vand{\oe}uvre-l\`es-Nancy, France.
                       Email: {\tt chassain@iecn.u-nancy.fr}}}

\begin{document}

\maketitle

\begin{center}

\vspace{3mm}

\rm
\end{center}

\begin{abstract}
We derive the distribution 
of the center of mass $S$ of the
integrated superBrownian excursion (ISE)
{from} the asymptotic distribution of the Wiener 
index for simple trees.
Equivalently, this is the distribution 
of the integral of a Brownian snake.
A recursion formula for the moments and asymptotics
for moments and tail probabilities are derived.

\vspace{3mm}

\vspace{3mm}
\textit{A.M.S.Classification.} 60K35  (primary), 
60J85 (secondary).
\end{abstract}

\section{Introduction}
The ISE (integrated superBrownian excursion)  is
 a random probability measure on  $\mathbb R^d$. 
The ISE was introduced by
David Aldous \cite{AldousISE} as an universal limit object 
 for random 
distributions of mass in $\mathbb R^d$: for instance,
Derbez \& Slade  \cite{DerbezSlade} 
proved that the ISE is the limit of 
lattice trees for $d\ge 8$. 
The ISE can be seen as the  limit of a suitably  renormalized
spatial branching process (cf.\ \cite{CnS,Mn'M}), or
 equivalently, as an embedding of the
continuum random  tree (CRT)  in $\mathbb R^d$. The ISE is  surveyed 
in \cite{Slade}.

Formally, the ISE 
is a random
variable, denoted $\ISE$, with value in the set of
 probability measures on
$\mathbb R^d$. 
In Section \ref{BrownianSnake}, 
we  give  a precise description
 of $\ISE$ in terms of the \textit{Brownian snake}, following 
 \cite[Ch. IV.6]{LEG}.
As noted in \cite{AldousISE}, even in the case $d=1$, 
where the support of $\ISE$ is almost surely
 a (random) bounded interval
denoted $[R,L]$,
 little was known about
the  distributional properties of 
 elementary statistics
 of $\ISE$, such as
$R$, $L$, or the center of mass 
\[S=\int x\ISE(dx).\]
This is still true, but recently $R-L$ was shown  \cite{CnS}
 to be the law 
of the radius of some model of random geometries
called \textit{fluid lattices} in
Quantum geometry  \cite{ambjorn}, or
 \textit{random quadrangulations} in combinatorics,
and the law of $R$ and $L$ were investigated by 
J.F.  Delmas  \cite{Delmas}. The joint law of 
$(R,L)$ is still unknown, as far as we know.

The Wiener index $w(G)$  of a connected
graph $G$ with set of vertices $V(G)$ is defined as
\[w(G)=\sum_{(u,v)\in V(G)}d(u,v),\]
in which $d(u,v)$ denotes the distance between $u$ and $v$
in the graph. For $G$ a random simple tree  with $n$ nodes,
$w(G)$, suitably normalized, is asymptotically distributed as
$\xi -\eta$, in which $\xi$ and $\eta$ are the following
 statistics
of the normalized Brownian excursion $\pa{e(t)}_{0\le t\le 1}$:
\begin{align}
\xi &= 2\int_0^1\,e(t)\,dt \label{xi}
\\
\eta &= 4\int_{0\le s<t\le 1}\min_{s\le u\le t}e(u)\,ds\,dt
\label{eta}
\end{align}
(cf.\ \cite{SJWiener}, where the joint moments of $(\eta,\xi)$
 are computed
{from} explicit formulas holding for random binary trees).

In this paper, we derive asymptotics, and
explicit induction formulas, for the moments of $S$,
{from} similar results for $\eta$.

\section{The ISE and the Brownian snake}
\label{BrownianSnake}

We recall briefly the description of $\ISE$ in terms of
the Brownian snake, from \cite[Ch. IV.6]{LEG}
 (see also \cite{BCHS,DZ,SERLET}).
The \textit{lifetime process} $\zeta=\pa{\zeta(t)}_{t\in T}$
of the Brownian snake
is a stochastic process with values in $[0,+\infty)$. 
Let the  Brownian snake with
lifetime $\zeta$ be denoted
\[W=\pa{W_s(t)}_{0\le s\le 1,\ 0\le t\le \zeta(s)}.\]
The Brownian snake can be seen
as a description of a "continuous" population $T$, through its genealogical tree
and the positions of its members. The lifetime $\zeta$ specifically describes
the genealogical tree, and $W$ describes the spatial motions of 
the members of population $T$. 

A member of the population is encoded by the time $t$ it is  visited 
by the contour traversal of the genealogical tree, $\zeta(t)$ being
the height of member 
 $t\in T$
in the genealogical tree ($\zeta(t)$ can be seen as 
the "generation" $t$ belongs to, or the time when $t$ is living).  Let 
\[C(s,t)=\min_{s\le u\le t} \zeta(u),\qquad
s\wedge t=\,\build{\textrm{argmin}}{s\le u\le t}{} \zeta(u).\]
Due to the properties of the contour traversal of a tree,
any element of $s\wedge t$ is a label for the more recent ancestor common to 
$s$ and $t$, and the distance between $s$ and $t$ in the genealogical tree
is 
\[d(s,t)=\zeta(s)+\zeta(t)-2C(s,t).\] 
If it is not a leaf of the tree, a member of the population is visited several
times ($k+1$ times if it has $k$ sons), so it has several labels:
$s$ and $t$ are two labels of the same member of the population
if $d(s,t)=0$, or equivalently if $s\wedge t\supset\{s,t\}$. 
Finally, $s$ is an ancestor of 
$t$ iff $s\in s\wedge t$. In this interpretation,
$W_s(u)$ is the position of the ancestor of $s$ living at time $u$, and 
\[
\widehat W_s=W_s\pa{\zeta(s)},  \qquad
 s\in T,
\]
is the position of $s$. Before time $m=C(s_1,s_2)$,
 $s_1$ and $s_2$ share the same ancestor,
entailing that 
\begin{equation}
\label{compatible}
\pa{W_{s_1}(t)}_{0\le t\le m}=
\pa{W_{s_2}(t)}_{0\le t\le m}.
\end{equation}
Obviously there is some redundancy in this description:
 it turns out that the full Brownian snake
can be recovered {from}  the pair
$\pa{\widehat W_s,\zeta(s)}_{0\leq s\leq1}$ (see \cite{Mn'M} for a
complete discussion of this).

In the general setting  \cite[Ch. IV]{LEG}, the spatial motion of
 a member of the population is any Markov proces with cadlag paths.
In the special case of the ISE, this spatial motion is a $d$-dimensional
Brownian motion:
\begin{itemize}
\item[a)]
for all $0\leq s\leq 1$, $t\rightarrow W_s(t)$ is a standard linear
 Brownian
motion started at 0,  defined for $0\le t\le \zeta(s)$
;
\item[b)] conditionally, given $\zeta$, the application $s\rightarrow W_s(.)$ is a
path-valued Markov process with transition function
defined as follows:
  for $s_1<s_2$, 
conditionally given $W_{s_1}(.)$,
 $\pa{W_{s_2}(m+t)}_{0\le t\le
\zeta(s_2)-m}$ is a standard Brownian motion starting {from} $W_{s_1}(m)$,
independent of $W_{s_1}(.)$.
\end{itemize}

The lifetime $\zeta$ is usually a reflected linear Brownian motion \cite{LEG},
defined on $T=[0,+\infty)$. However,
  in the case  of the ISE,
\[\zeta=2e,\] 
in which $e$ denotes
the normalized Brownian excursion, 
or 3-dimensional Brownian bridge, defined on $T=[0,1]$. 
With this choice of  $\zeta$,
the genealogical tree is the CRT (see \cite{AldousISE}), 
and the Brownian snake can be seen as 
an embedding of the CRT 
in $\bbR^d$.  
We can now give the definition of the ISE
in terms of the Brownian snake with lifetime $2e$ \cite[Ch. IV.6]{LEG}:
\begin{defn} The ISE
$\ISE$ is the occupation measure of $\widehat W$.
\end{defn} 

Recall that the occupation measure $\ISE$ of a process $\widehat W$ is 
defined by the relation:
\begin{equation}\label{occu}
  \int_{\mathbb R}f(x)\ISE(dx)=\int_0^1f\pa{\widehat W_s}ds,
\end{equation}
holding for any measurable test function $f$.

\begin{rem}
It might seem more natural to consider the Brownian snake with lifetime $e$
instead of $2e$, but we follow the normalization in Aldous  \cite{AldousISE}.
If we used the Brownian snake with lifetime $e$ instead,
$W$, $\ISE$ and $S$ would be scaled by $1/\sqrt2$.
\end{rem}

\section{The basic identity}
\begin{thm}
\label{idloi}
 Let $N$ be a standard Gaussian random
 variable, independent of $\eta$. Then
\[S\loi \sqrt{\eta}\  N.\]
\end{thm}
\begin{proof}
Based on the short account of the Brownian snake theory in
\refS{BrownianSnake}, 
the proof is now pretty easy. Specializing \eqref{occu} to $f(x)=x$, 
we obtain a
representation of $S$:
\[S=\int_0^1\widehat W_s\ ds,\]
which is the starting point of our proof.
We have also,   directly {from} the definition of
 the Brownian snake,

\begin{prop}
\label{gaussianprocess} Conditionally, given $e$,
$\pa{\widehat W_s}_{0\le s\le 1}$  is a Gaussian process
whose  covariance is  $C(s,t)=2\min_{s\le u\le t} e(u)$, $s\le t$. 
\end{prop}
\begin{proof} With the notation $\zeta=2e$ and 
$m=C(s_1,s_2)=2\min_{s_1\le u\le s_2} e(u)$, we have, 
conditionally, given $e$, for $s_1\le s_2$,
\begin{eqnarray*}
\cov{\widehat W_{s_1},\widehat W_{s_2}}
 &=& \cov{W_{s_1}(\zeta(s_1)),
W_{s_2}(\zeta(s_2))-W_{s_2}(m)+W_{s_2}(m)}
\\
&=& \cov{W_{s_1}(\zeta(s_1)),W_{s_2}(m)}
\\
&=& \cov{W_{s_1}(\zeta(s_1)),W_{s_1}(m)}
\\
&=& m,
\end{eqnarray*}
in which b) yields the second equality,
(\ref{compatible}) yields the third one, and a) yields the fourth equality.
\end{proof}

As a consequence of Proposition \ref{gaussianprocess},
conditionally given $e$, $S$ is centered Gaussian with variance
\[\int_{[0,1]^2}C(s,t)\ ds\,dt=\eta.\]
This last statement is equivalent to Theorem \ref{idloi}.
\end{proof}

\begin{rem} The $d$-dimensional analog of
Theorem \ref{idloi} is also true, by the same proof.
\end{rem}

\section{The moments}
We know that the odd moments of $S$ vanish.
\begin{thm}
\label{moments}
For $k\ge 0$,
\begin{equation}
\label{esk}
\esp{S^{2k}} = \frac{(2k)!\sqrt{\pi}}{2^{(9k-4)/2}\Gamma\pa{(5k-1)/2}}\; a_k,
\end{equation}
in which  $a_k$ is defined by  $a_1=1$, and, for $k\ge 2$,
\begin{equation}
\label{induction1}
a_k=2(5k-4)(5k-6)a_{k-1}+\sum_{i=1}^{k-1}a_ia_{k-i}.
\end{equation}
\end{thm}
\begin{proof}
{}{From} Theorem \ref{idloi},  we derive at once that
\[\esp{ S^{2k}}=\mathbb E\cro{\eta^{k}}\mathbb E\cro{N^{2k}}.\]
 As a special case of \cite[Theorem 3.3]{SJWiener}
(where $a_k$ is denoted $\omega^*_{0k}$),
\begin{equation}
\label{etak}
\esp{\eta^{k}} = \frac{k!\sqrt{\pi}}{2^{(7k-4)/2}\Gamma\pa{(5k-1)/2}}\;  a_k,
\end{equation}
and the result follows since $\esp{N^{2k}}=(2k)!/(2^kk!)$.
\end{proof}

In particular, see again
\cite[Theorem 3.3]{SJWiener},
\begin{align*}
 \esp{S^2}&=\esp{\eta}=\sqrt{\pi/8},\\ 
\esp{S^4}&=3\esp{\eta^2}=7/5,
\end{align*}
as computed by Aldous \cite{AldousISE} by a different method.
(Aldous' method extends to higher moments too, but the calculations 
quickly become complicated.)

We have the following asymptotics for the moments.
\begin{thm}
\label{asympt}
For some constant $\beta=0.981038\dots$ 
we have, as $k\to\infty$,
\begin{align}
a_k&\sim \beta\, 50^{k-1}\pa{k-1}!^2,  \label{2a}
\\
\esp{\eta^k} &\sim \frac{2\pi^{3/2}\beta}{5} k^{1/2}
\pa{5e}^{-k/2} k^{k/2} \label{2n},
\\
\esp{S^{2k}} &\sim \frac{2\pi^{3/2}\beta}{5} (2k)^{1/2}
\pa{10e^{3}}^{-2k/4} (2k)^{\frac34\cdot 2k} \label{2s}.
\end{align}
\end{thm}
\begin{proof} 
Set
\begin{equation*}
b_k=\frac{a_k}{\ 50^{k-1}\pa{k-1}!^2}.
\end{equation*}
We have $b_1=a_1=1$ and, {from} \eqref{induction1},
$b_2=b_3=49/50=0.98$.
For $k\ge 3$, (\ref{induction1}) translates to
\begin{eqnarray*}
b_k&=&\pa{1-\frac 1{25(k-1)^2}}b_{k-1}+
\frac{\sum_{i=1}^{k-1}a_ia_{k-i}}{50^{k-1}\pa{k-1}!^2}\\
&=&b_{k-1}+
\frac{\sum_{i=2}^{k-2}a_ia_{k-i}}{50^{k-1}\pa{k-1}!^2},
\end{eqnarray*}
Thus $b_k$
 increases for $k\ge 3$. 
We will show that $b_k<1$ for all $k>1$; thus 
$\beta=\lim_{k\to\infty} b_k$
exists and \eqref{2a} follows.
Stirling's formula and \eqref{etak}, \eqref{esk} then  yield
\eqref{2n} and \eqref{2s}.

More precisely,
we show by induction that for $k\ge3$,
\[b_k\le 1-\frac1{25(k-1)}.\]
This holds for $k=3$ and $k=4$. For $k\ge5$ we have by the induction  
assumption $b_j\le1$ for $1\le j <k$ and thus
\begin{equation}
\label{sk}
  \begin{split}
s_k&=\frac{\sum_{i=2}^{k-2}a_ia_{k-i}}{50^{k-1}\pa{k-1}!^2}\\
&=\frac1{50}\sum_{i=2}^{k-2}\frac{\pa{i-1}!^2\pa{k-i-1}!^2}{\pa{k-1}!^2}
b_i b_{k-i}\\
&\le\frac1{50\pa{k-1}^2\pa{k-2}^2}\ \ \pa{2+
(k-5)\frac{4}{\pa{k-3}^2}}\\
&\le\frac1{25\pa{k-1}\pa{k-2}^2\pa{k-3}}.
  \end{split}
\end{equation}
The induction  follows. Moreover, it follows easily {from}
\eqref{sk} that $b_k < \beta < b_k+\frac1{75}(k-2)^{-3}$.

To obtain the numerical value of $\beta$, we write, somewhat more sharply,
where $0\le \theta \le 1$,
\begin{multline*}
s_k
= \frac1{25\pa{k-1}^2\pa{k-2}^2} b_2 b_{k-2}
+ \frac4{25\pa{k-1}^2\pa{k-2}^2(k-3)^2} b_3 b_{k-3}\\
+\theta (k-7)\frac{36}{25\pa{k-1}^2\pa{k-2}^2(k-3)^2(k-4)^2}
\end{multline*}
and sum over $k>n$ for $n=10$, say, using
$b_{n-1}<b_{k-2}<\beta$ and $b_{n-2}<b_{k-3}<\beta$ for $k>n$.
It follows (with \texttt{Maple}) by this and exact computation of
$b_1,\dots,b_{10}$
that $0.981038<\beta <0.9810385$; we omit the details.
\end{proof}

\begin{rem}
For comparison, we give the corresponding result for
$\xi$ defined in \eqref{xi}. 
There is a simple relation, 
discovered by
Spencer \cite{Spencer} and Aldous \cite{Aldous:mult},
between its moments 
and Wright's constants in the enumeration of connected graphs with $n$
vertices and $n+k$ edges \cite{Wright},
and the well-known asymptotics of the latter lead, see
\cite[Theorem 3.3 and (3.8)]{SJWiener}, 
to 
\begin{equation}
\label{xix}
\esp{\xi^k} \sim 3\sqrt{2}{k}(3e)^{-k/2} k^{k/2}.
\end{equation}
\end{rem}

\section{Moment generating functions and tail estimates}

The moment asymptotics yield asymptotics for the moment generating
function
$\esp{e^{t S}}$ and the tail probabilities $\pr{S>t}$ as $t\to\infty$.
For completeness and comparison, we also include corresponding results
for $\eta$.
We begin with a standard estimate.
\begin{lem} \label{L1}
\begin{romenumerate}
\item
If $\gamma>0$ and $b\in\bbR$, then, as $x\to\infty$,
\begin{equation*}
  \sum_{k=1}^\infty k^b k^{-\gamma k} x^k
\sim
\parfrac{2\pi}{\gamma}^{1/2} \bigpar{e^{-\gamma}x}^{(b+1/2)/\gamma}
e^{\gamma \smallpar{e^{-\gamma}x}^{1/\gamma}}
\end{equation*}
\item
If $-\infty<\gamma<1$ and $b\in\bbR$, then, as $x\to\infty$,
\begin{equation*}
  \sum_{k=1}^\infty \frac{k^b k^{\gamma k} x^k}{k!}
\sim
(1-\gamma)^{-1/2} \bigpar{e^{\gamma}x}^{b/(1-\gamma)}
e^{(1-\gamma) \smallpar{e^{\gamma}x}^{1/(1-\gamma)}}
\end{equation*}
\end{romenumerate}
The sums over even $k$ only are asymptotic to half the full sums.
\end{lem}
\begin{proof}[Sketch of proof]
(i).
This is standard, but since we have not found a precise
 reference, we sketch the argument.
Write 
$ k^b k^{-\gamma k} x^k=e^{f(k)}=k^b e^{g(k)}$
where $g(y)=-\gamma y \ln y + y \ln x$ and $f(y)=g(y)+b\ln y$.
The function $g$ is concave with a maximum at
$y_0=y_0(x)=e^{-1}x^{1/\gamma}$.
A Taylor expansion yields

\begin{align*}
y_0^{-1/2}e^{-f(y_0)}\sum_{k=1}^\infty k^b k^{-\gamma k} x^k
&=
y_0^{-1/2}\int_0^\infty e^{f(\ceil{y})-f(y_0)}dy\\
&=\int_{-y_0^{1/2}}^\infty e^{f(\ceil{y_0+sy_0^{1/2}})-f(y_0)}ds
\\
&\to
\int_{-\infty}^\infty e^{-\gamma s^2/2}ds =\parfrac{2\pi}\gamma^{1/2}.
\end{align*}

(ii). Follows by (i) and Stirling's formula.
\end{proof}

Combining \refT{asympt} and \refL{L1}, 
we find the following asymptotics for
$\esp{e^{t\eta}} =\sum_k \esp{\eta^k}t^k/k!$
and
$\esp{e^{tS}} =\sum_{k \text{ even}} \esp{S^k}t^k/k!$.
\begin{theorem}
  \label{Tmgf}
As $t\to\infty$,
\begin{align}
\esp{e^{t\eta}} 
&\sim 
\frac{(2\pi)^{3/2}\beta}{5^{3/2}} t
e^{t^2/10}, 
\\
\esp{e^{tS}} &\sim \frac{2^{1/2}\pi^{3/2}\beta}{5^{3/2}}t^2 
e^{t^4/40}.
\end{align}
\end{theorem}
\begin{proof}
  For $\eta$ we take $b=1/2$, $\gamma=1/2$, $x=(5e)^{-1/2}t$
in \refL{L1}(ii);
for $S$ we take
$b=1/2$, $\gamma=3/4$, $x=(10e^3)^{-1/4}t$.
\end{proof}

The standard argument with Markov's inequality yields upper bounds for
the tail probabilities {from} \refT{asympt} or~\ref{Tmgf}.
\begin{theorem}
  \label{Ttail1}
For some constants $K_1$ and $K_2$ and all $x\ge1$, say,
\begin{align}
\pr{\eta>x} 
&\le K_1 x\exp\bigpar{-\tfrac52 x^2},
\label{etax}
\\
\pr{S>x} 
&\le K_2 x^{2/3} \exp\bigpar{-\tfrac34 10^{1/3} x^{4/3}}.
\label{sx}
\end{align}
\end{theorem}
\begin{proof}
  For any even $k$ and $x>0$, $\pr{|S|> x} \le x^{-k}\esp{S^k}$.
We use \eqref{2s} and 
optimize the resulting exponent by
choosing $k=10^{1/3}x^{4/3}$, rounded to an even integer.
This yields \eqref{sx}; we omit the details.
\eqref{etax} is obtained similarly {from} \eqref{2n}, using
$k=\floor{5x^2}$. 
\end{proof}

\begin{rem}
  The proof of \refT{Ttail1} shows that any 
$K_1>\xfrac{2\pi^{3/2}\beta}{5^{1/2}}\approx 4.9$
and
$K_2>\xfrac{10^{1/6}\beta\pi^{3/2}}{5}\approx 1.6$
will do for large $x$. Alternatively, we could use \refT{Tmgf} and
$\pr{S>x}<e^{-tx}\esp{e^{tS}}$ for $t>0$ and so on; this would yield
another proof of \refT{Ttail1} with somewhat inferior values of $K_1$
and $K_2$. 
\end{rem}

The bounds obtained in \refT{Ttail1} are sharp up to factors $1+o(1)$
in the exponent, as is usually the case for estimates derived by
this method. 
For convenience we state a general theorem,
reformulating results
by Davies \cite{Davies} and Kasahara \cite{Kasahara}.

\begin{theorem} \label{Tabc}
  Let $X$ be a random variable, let $p>0$ and let $a$ and $b$ be
  positive real numbers related by 
$a=1/(p e b^{p})$ or, equivalently,
$b=(pea)^{-1/p}$.
\begin{romenumerate}
\item
If $X\ge0$ a.s., then
\begin{equation}
  \label{ax}
-\ln\P(X>x) \sim ax^{p} \qquad \text{as $x\to\infty$}
\end{equation}
is equivalent to
\begin{equation}
  \label{bx}
\bigpar{\E X^r}^{1/r} \sim b r^{1/p} \qquad \text{as $r\to\infty$}.
\end{equation}
Here $r$ runs through all positive reals; equivalently, we can
restrict $r$ in \eqref{bx} to integers or even integers.
\item
If $X$ is a symmetric random variable, then \eqref{ax} and \eqref{bx}
are equivalent, 
where $r$ in \eqref{bx} runs through 
even integers.
\item
If $p>1$, then, for any $X$, \eqref{ax} is equivalent to 
\begin{equation}
  \label{cx}
\ln\bigpar{\E e^{tX}} \sim c t^{q} \qquad \text{as $t\to\infty$},
\end{equation}
where $1/p+1/q=1$ and 
$c=q^{-1}(pa)^{-(q-1)}
=q^{-1}e^{q-1}b^{q}$.
(This can also be written in the symmetric form 
$(pa)^q(qc)^p=1$ and as $b=e^{-1/p}(qc)^{1/q}$.)
Hence, if $X\ge0$ a.s., or if $X$ is symmetric and $r$ restricted to
even integers, \eqref{cx} is also equivalent to \eqref{bx}.
\end{romenumerate}
\end{theorem}

\begin{proof}
For $X\ge0$,
(i) and (iii) 
are immediate special cases of more general results by
Kasahara \cite[Theorem 4 and Theorem 2, Corollary 1]{Kasahara}, see 
also \cite[Theorem 4.12.7]{regular};
the difficult implications \eqref{bx}$\implies$\eqref{ax}
and \eqref{cx}$\implies$\eqref{ax} were earlier proved by Davies
\cite{Davies} (for $p>1$, which implies the general case of (i) by
considering a power of $X$). 
Moreover, \eqref{cx}$\implies$\eqref{ax} follows also {from} the
G\"artner--Ellis theorem
\cite[Theorem 2.3.6]{DemboZ} applied to $n^{-1/p}X$.
Note that, assuming $X\ge0$,
$(\E X^r)^{1/r}$ is increasing in $r>0$,
which implies that \eqref{bx} for (even) integers is equivalent to
\eqref{bx} for all real $r$.

(ii) follows {from} (i) applied to $|X|$, and (iii) for general $X$
follows by considering $\max(X,0)$.
\end{proof}

We thus obtain {from} Theorem \ref{asympt} or \ref{Tmgf}
the following estimates, less precise than
\refT{Ttail1} but including both upper and lower bounds. 
\begin{thm}
  \label{Ttail2}
As $x\to\infty$,
\begin{align}
 \ln\bigpar{\pr{\eta>x}}& \sim -\frac52 x^2 , \label{eta2}
\\
 \ln\bigpar{\pr{S>x}} & 
\sim -\frac 34 10^{1/3}\, x^{4/3}. \label{S2}
\end{align}
\end{thm}
\begin{proof}
For $\eta$ we use \refT{Tabc} with $p=q=2$, $b=(5e)^{-1/2}$,
$a=5/2$ and $c=1/10$.
For $S$ we have $p=4/3$, $q=4$, $b=(10e^3)^{-1/4}$,
$a=3\cdot10^{1/3}/4$ and $c=1/40$.
\end{proof}

\begin{rem}
\eqref{eta2} can also be proved using the representation \eqref{eta}
and
large deviation theory for Brownian excursions,
cf. \cite[\S 5.2]{DemboZ} and \cite{DZ}.
The details may perhaps appear elsewhere.
\end{rem}

\begin{rem} 
\eqref{S2} can be compared to the tail estimates
in \cite{AldousISE} for the density function of 
$\widehat W _U$, the 
value of the Brownian snake 
at a random point $U$,
which in particular gives
$$
 \ln\bigpar{\pr{\widehat W_U>x}} 
 \sim\  - 3\cdot 2^{-5/3}\, x^{4/3}.
$$
\end{rem}

\begin{rem}
For $\xi$ we can by \eqref{xix} use \refT{Tabc} with $p=q=2$, $b=(3e)^{-1/2}$,
$a=3/2$ and $c=1/6$. In \cite{SJWiener}, the variable of main interest
is neither $\xi$ nor $\eta$ but $\zeta=\xi-\eta$. By Minkowski's
inequality
${\esp{\xi^k}}^{1/k} 
\le {\esp{\zeta^k}}^{1/k}  + {\esp{\eta^k}}^{1/k} 
$
and \eqref{2n}, \eqref{xix} follows 
\begin{equation*}
\frac1{\sqrt{3e}} - \frac1{\sqrt{5e}} 
\le   \liminf_{k\to\infty}  \frac{\bigpar{\esp{\zeta^k}}^{1/k} }{k^{1/2}}
\le   \limsup_{k\to\infty}  \frac{\bigpar{\esp{\zeta^k}}^{1/k} }{k^{1/2}}
\le \frac1{\sqrt{3e}}. 
\end{equation*}
This leads to 
asymptotic upper and lower bounds for $\ln(\pr{\zeta>x})$ too by
\cite{Davies} or \cite{Kasahara}.
We can show that 
$\lim_{k\to\infty} k^{-1/2} \bigpar{\esp{\zeta^k}}^{1/k} $ and 
$\lim_{x\to\infty} x^{-2}\ln(\pr{\zeta>x})$ 
exist, but do not
know their value.
\end{rem}

\section{Concluding remarks}
The center of mass  of the ISE 
 turns out to be related to the Wiener index of simple trees,
but note that  $S$ is also related
to the Wiener index of random planar quadrangulations:
let $(Q_n,(b_n,e_n))$  denote the uniform choice of a
 quadrangulation
 with $n$ faces, and of a "marked" oriented edge in it, and set
\[
\mathcal W_n=\sum_{x\in Q_n} d(x,b_n).
\]
As a consequence of  \cite{CnS},
$n^{-5/4}\mathcal W_n$ converges weakly to
$c\cdot(S-L)$, where $L$ is the left endpoint
of the support of $\ISE$, and $c$ is a known constant.
The joint law of $(S,L)$ is not known, as far as we know.

\begin{ack}
  We thank Jim Fill for helpful comments.
\end{ack}

\newcommand\vol{\textbf}
\newcommand\jour{\emph}
\newcommand\book{\emph}
\newcommand\inbook{In \emph}


\begin{thebibliography}{99}


\bibitem{AldousISE} D. Aldous.
\newblock Tree-based models for random distribution of mass.
\newblock {\em Journal of Statistical Physics} \vol{73} (1993), 625-641.

\bibitem{Aldous:mult} D. Aldous.
Brownian excursions, critical random graphs
and the multiplicative coalescent.
\emph{Ann. Probab.} \vol{25} (1997), no. 2, 812--854.


\bibitem{ambjorn} J. Ambj\o rn, B. Durhuus and T. J\'onsson.
\newblock {\em Quantum gravity, a statistical field theory approach.}
\newblock Cambridge Monographs on Mathematical Physics,  1997.

\bibitem{regular} N. H. Bingham, C. M. Goldie and
 J. L. Teugels.
\newblock {\em Regular variation.}  First edition.
\newblock Encyclopedia of Mathematics and its Applications, 27. 
 Cambridge University Press, Cambridge, 1987.


\bibitem{BCHS} C.
Borgs, J. Chayes, R. van der Hofstad and G. Slade.
\newblock
Mean-field lattice trees.
        On combinatorics and statistical
mechanics.
\newblock {\em Ann. Comb.} \vol3 (1999), no. 2--4, 205--221.


\bibitem{CnS}
P. Chassaing and G. Schaeffer.
\newblock Random Planar Lattices and
Integrated SuperBrownian Excursion.
\newblock To appear in
\jour{Probab. Theory Related Fields}, 
2002.

\bibitem{Davies}
L. Davies. 
Tail probabilities for positive random variables with entire
characteristic functions of very regular growth. 
\jour{Z. Angew. Math. Mech.} \vol{56} (1976), no. 3, T334--T336.

\bibitem{Delmas}
J. F. Delmas.
\newblock Computation of moments for the length of the one dimensional ISE
support. \newblock {\em Preprint}, 2002, available at 
\newline
\texttt{ http://cermics.enpc.fr/\string~delmas/publications.html}.

\bibitem{DemboZ}
A. Dembo and O. Zeitouni. 
\book{Large deviations techniques and applications.} 
Jones and Bartlett Publishers, Boston, MA, 1993. 



\bibitem{DZ}  
A. Dembo and O. Zeitouni.
\newblock Large  deviations for random
distribution of mass.
\newblock {\em Random discrete structures
(Minneapolis, MN, 1993)},  45--53,
IMA Vol. Math. Appl.,
76,
Springer, New York, 1996.


\bibitem{DerbezSlade} E. Derbez and G. Slade.
\newblock The scaling limit of lattice trees
in high dimensions.
\newblock \jour{Commun. Math. Phys.} \vol{193} (1998), 69--104.


\bibitem{SJWiener}
S. Janson.
\newblock The Wiener index of simply generated random trees.
\emph{Random Struct. Alg.} \vol{22} (2003), no. 4, 337--358. 


\bibitem{Kasahara}
Y. Kasahara. 
Tauberian theorems of exponential type.  
\jour{J. Math. Kyoto Univ.}  \vol{18}  (1978), no. 2, 209--219

\bibitem{LEG}  J.-F. Le Gall.  \newblock {\em
Spatial branching processes,
      random snakes and partial differential
equations.}
Lectures in Mathematics ETH Z\"urich, Birkh\"auser
Verlag, Basel, 1999.


\bibitem{Mn'M}
J.-F. Marckert and A. Mokkadem.
\newblock States spaces of the snake and of its tour -- Convergence of
the discrete snake.
\newblock To appear in \jour{J. Theoret. Probab.}, 2002.




\bibitem{SERLET} L. Serlet.
\newblock A large
deviation principle for the Brownian snake.
\newblock {\em Stochastic
Process. Appl.} \vol{67} (1997), no. 1, 101--115.


\bibitem{Slade} G. Slade.
\newblock Scaling limits and
Super-Brownian motion.
\newblock {\em Notices of the AMS}  \vol{49} (2002), no. 9, 1056--1067.

\bibitem{Spencer} 
J. Spencer.
Enumerating graphs and Brownian motion.
\emph{Commun. Pure Appl. Math.} \vol{50} (1997), 
 291--294.

\bibitem{Wright}
E.M. Wright.
The number of connected sparsely edged graphs. 
\emph{J. Graph Th.} \vol{1} (1977),
317--330.

\end{thebibliography}

\end{document}